\documentclass[10pt]{article}
\usepackage{amsmath,amssymb}
\usepackage{enumerate}
\usepackage{fullpage}

\allowdisplaybreaks[1]

\newtheorem{theorem}{Theorem}[section]

\newtheorem{lemma}[theorem]{Lemma}
\newtheorem{proposition}[theorem]{Proposition}
\newtheorem{corollary}[theorem]{Corollary}
\newtheorem{definition}[theorem]{Definition}

\newtheorem{example}[theorem]{Example}

\newcommand{\End}{\operatorname{End}}

\newcommand{\Aut}{\operatorname{Aut}}

\newcommand{\Map}{\operatorname{Map}}

\newenvironment{proof}{\par\noindent{\bf Proof.}}{$\qed$\par\bigskip}
\newcommand{\qed}{\enspace\vrule  height6pt  width4pt  depth2pt}
\usepackage{color}

\begin{document}

\title{Corrigendum and Addendum to ``Structure monoids of set-theoretic solutions of the Yang--Baxter equation''
\thanks{The first author was partially supported by the grant MINECO PID2020-113047GB-I00 (Spain).
The
second author was supported in part by Onderzoeksraad of Vrije
Universiteit Brussel and Fonds voor Wetenschappelijk Onderzoek
(Belgium). The third author is supported by Fonds voor Wetenschappelijk Onderzoek (Flanders), via an FWO Aspirant-mandate.}
}\author{Ferran Ced\'o \and Eric Jespers \and Charlotte Verwimp}
\date{}

\maketitle

\begin{abstract}
One of the results in our article, which appeared in Publ. Mat. 65 (2021), 499--528, is that the structure monoid $M(X,r)$ of a left non-degenerate solution $(X,r)$ of the Yang-Baxter Equation  is a left semi-truss, in the sense of Brzezi\'nski, with an additive structure monoid that is close to being a normal semigroup. Let $\eta$ denote the least left cancellative congruence on the additive monoid $M(X,r)$.
It is then  shown  that 
$\eta$ also is a congruence on the multiplicative monoid $M(X,r)$ and that 
the left cancellative epimorphic image $\bar{M}=M(X,r)/\eta$ 
inherits a semi-truss structure and thus one obtains a natural  left non-degenerate solution of the Yang-Baxter equation on $\bar{M}$. 
Moreover, it restricts to the original solution $r$ for some interesting classes, in particular if $(X, r)$ is irretractable. 
The proof contains a gap. In the first part of the paper we correct this mistake by introducing a new left cancellative congruence $\mu$ on the additive monoid $M(X,r)$ and show that it also yields a  left cancellative congruence on the multiplicative monoid $M(X,r)$ and we obtain a semi-truss structure on $M(X,r)/\mu$ that also yields a natural left non-degenerate solution.

In the second part of the paper we start from the least left cancellative congruence $\nu$ on the multiplicative monoid $M(X,r)$ and show that it also is a congruence on the additive monoid $M(X,r)$ in case $r$ is bijective. If, furthermore, $r$ is left and right non-degenerate and bijective then $\nu =\eta$, the least left cancellative congruence on the additive monoid $M(X,r)$, extending an earlier result of Jespers, Kubat and Van Antwerpen to the infinite case.
\end{abstract}

\bigskip

\bigskip

\section{Introduction}
We have detected a mistake in the proof of \cite[Lemma 5.5]{CJV}. What is correctly proved is the following result for a left non-degenerate solution $(X,r)$ of the Yang-Baxter equation (YBE) with structure monoid $M=M(X,r)$.  Write $r(x,y)=(\sigma_x(y), \gamma_y(x))$.
Thus, all $\sigma_x$ are bijective maps.
Its additive structure is denoted by $(M,+)$ and its multiplicative structure  by $(M,\circ)$. The least cancellative congruence on $(M,+)$ is denoted  by $\eta$.
Let  $\lambda'\colon (M,\circ) \rightarrow \End (M,+)$ denote the unique monoid homomorphism such that $\lambda' (x)(y)=\sigma_x (y)$ for $x,y \in X$  (see Proposition~3.1 in \cite{CJV}).

\begin{lemma}\label{eta}
	With the same notation as in \cite[Lemma 5.5]{CJV} we have $\eta=\eta'$. Furthermore, for all $z\in M$,
	$$\eta=\{(\lambda'_z(a),\lambda'_z(b))\mid (a,b)\in \eta\}=\{((\lambda'_z)^{-1}(a),(\lambda'_z)^{-1}(b))\mid (a,b)\in \eta\}.$$  
\end{lemma}

We do not know whether $\lambda'_a=\lambda'_b$, for all $(a,b)\in \eta$ and whether $\eta$ is a congruence on $(M,\circ)$. As a consequence \cite[Remark 5.6, Corollary 5.9 and Corollary 5.10]{CJV} are not proved. Therefore \cite[Question 5.7]{CJV} and the definition of injective left non-degenerate solution of the YBE given in \cite{CJV} have no sense.
In Section 2 we will introduce a new congruence on $(M,+)$ and prove a correct version of the listed corollaries.

In Section 3,  we start from the least left cancellative congruence $\nu$ on the multiplicative  monoid $(M,\circ)$ and show that it also is a congruence on the additive  monoid $(M,+)$ in case $r$ is bijective. If furthermore $r$ is left and right non-degenerate and bijective, then $\nu =\eta$, the least left cancellative congruence on the additive  monoid $(M,+)$, extending an earlier result of Jespers, Kubat and Van Antwerpen to the infinite case.

\section{Correction of \cite[Section 5]{CJV}}
In this section, we shall introduce a new congruence $\mu$ on $(M,+)$ such that it also is a congruence on $(M,\circ)$, $(M,+)/\mu$ is left cancellative and $((\lambda'_a)^{\varepsilon}(b), (\lambda'_{a'})^{\varepsilon}(b'))\in \mu$, for all $(a,a'),(b,b')\in \mu$  and $\varepsilon\in\{ -1,1\}$. Furthermore, $\mu$ is the least binary relation on $M$ with these properties.

We first recall the definition of a left semi-truss.

\begin{definition} (Brzezi\'nski \cite{Br})
	A left semi-truss is a quadruple $(A,+,\circ,\phi)$ such that
	$(A,+)$ and $(A,\circ)$ are semigroups and $\phi\colon
	A\times A\longrightarrow A$ is a function such that
	$$a\circ (b+ c)= (a\circ b)+ \phi(a,c),$$
	for all $a,b,c\in A$.
\end{definition}

\begin{example}\cite[Example 5.2]{CJV}\label{ex1}
	{\rm Let $(X,r)$ be a left non-degenerate set-theoretic solution of
		the YBE (not necessarily bijective). Again write $r(x,y)=(\sigma_x(y),\gamma_y(x))$, for 
$x,y\in X$. 
		As stated in \cite[Section~3]{CJV},
		and with the same notation,  the map
		$$r'(x,y)=(y,\sigma_y\gamma_{\sigma^{-1}_x(y)}(x))$$ defines the left
		derived solution on $X$. Let $M=M(X,r)$ and $M'=A(X,r)=M(X,r')$ be
		the structure monoids of the solutions $(X,r)$ and $(X,r')$
		respectively. From \cite[Corollary~3.9 and
		Proposition~3.1]{CJV} we obtain a left action
		$\lambda'\colon (M,\circ )\longrightarrow \Aut(M',+)$ and a
		bijective $1$-cocycle $\pi\colon M\longrightarrow M'$ with respect
		to $\lambda'$ satisfying $\lambda'(x)(y)=\sigma_x(y)$ and
		$\pi(x)=x$, for all $x,y\in X$. We identify $M$ and $M'$ via $\pi$,
		that is, $a=\pi(a)$ for all $a\in M$. With this identification, we
		obtain  the operation $+$ on $M$, and $a\circ b=a+\lambda'_a(b)$,
		for all $a,b\in M$. Put $\phi(a,b)=\lambda'_a(b)$, for all $a,b\in
		M$. Then,
		$$a\circ (b+c)=a+\lambda'_a(b+c)=a+\lambda'_a(b)+\lambda'_a(c)=(a\circ b)+\phi(a,c),$$
		for all $a,b\in M$.
		Furthermore, $M + a \subseteq a+M$, for all $a\in M$. Hence
		$(M,+,\circ,\phi)$ is a left semi-truss. Note that if, furthermore,
		$r$ is bijective then it easily can be verified that $(X,r')$ is a
		right non-degenerate solution and thus  $M+a =a+M$ for all $a\in M$;
		that is, $(M,+)$ consists of normal elements. As shown in \cite{JKV},
		this property is fundamental in the study of the associated
		structure algebra $KM(X,r)$, where $K$ is a field.}
\end{example}

We will use the assumptions and notations as in Example \ref{ex1}.

Let
$$\mu_0=\{ (a,b)\in M^2\mid \exists\; c\in M \mbox{ such that } c+a=c+b\}.$$
Note that $\mu_0$ is a reflexive and symmetric binary relation on
$M$. Let $\mu_1$ be its transitive closure, that is
$$\mu_1=\{(a,b)\in M^2\mid \exists\;  a_1,\dots,a_n\in M \mbox{ such that }(a,a_1),(a_1,a_2),\dots ,(a_n,b)\in \mu_0\}.$$
Thus $\mu_1$ is an equivalence relation on $M$. Let
\begin{align*}\mu_2=&\{ ((\lambda'_{z})^{\varepsilon}(a\circ c),(\lambda'_{z})^{\varepsilon}(b\circ c))\in M^2\mid z,c\in M,\, \varepsilon\in \{ -1,1\}\mbox{ and }(a,b)\in \mu_1\} ,\\
	\mu_3=&\{(a,b)\in M^2\mid \exists\; a_1,\dots,a_n\in M \mbox{ such that }(a,a_1),(a_1,a_2),\dots ,(a_n,b)\in \mu_2\},\\
	\mu_4=&\{(c+a+d,c+b+d)\in M^2\mid c,d\in M \mbox{ and } (a,b)\in \mu_3\}\\
	&\cup \{ (a,b)\in M^2\mid \exists\; c\in M\mbox{ such that
	}(c+a,c+b)\in \mu_3\}, 
\end{align*}
and for every $m\geq 1$ we
define
\begin{align*}\mu_{4m+1}&=\{(a,b)\in M^2\mid \exists\; a_1,\dots,a_n\in M \mbox{ such that }(a,a_1),(a_1,a_2),\dots ,(a_n,b)\in \mu_{4m}\}\\
	\mu_{4m+2}=&\{ ((\lambda'_{z})^{\varepsilon}(a\circ c),(\lambda'_{z})^{\varepsilon}(b\circ c))\in M^2\mid z,c\in M,\, \varepsilon\in\{ -1,1\}\mbox{ and }(a,b)\in \mu_{4m+1}\} ,\\
	\mu_{4m+3}=&\{(a,b)\in M^2\mid \exists\; a_1,\dots,a_n\in M \mbox{ such that }(a,a_1),(a_1,a_2),\dots ,(a_n,b)\in \mu_{4m+2}\},\\
	\mu_{4(m+1)}=&\{(c+a+d,c+b+d)\in M^2\mid c,d\in M \mbox{ and } (a,b)\in \mu_{4m+3}\}\\
	&\cup \{ (a,b)\in M^2\mid \exists\; c\in M\mbox{ such that
	}(c+a,c+b)\in \mu_{4m+3}\}.
\end{align*}

Note that
$\mu_n\subseteq\mu_{n+1}$, for all $n\geq 0$. Let
$\mu=\cup_{n=0}^{\infty}\mu_n$.

\begin{lemma}\label{congruence}
	With the above notation, we have that $\mu$ is a congruence on $(M,+)$ and  it also is a congruence on $(M,\circ)$. Furthermore, $(M,+)/\mu$ and $(M,\circ)/\mu$ are left cancellative monoids, and  
	$$( \lambda'_c(a),\lambda'_d(b)), ( (\lambda'_c)^{-1}(a),(\lambda'_d)^{-1}(b))\in\mu,$$
	for all $(a,b),(c,d)\in \mu$. 
\end{lemma}
\begin{proof}
	First we shall prove that $\mu$ is a congruence on $(M,+)$.
	Clearly $\mu$ is reflexive and symmetric because so is each
	$\mu_n$. Let $a,b,c\in M$ be  such that $(a,b),(b,c)\in \mu$. There
	exists a positive integer $m$ such that $(a,b),(b,c)\in
	\mu_{2m}$. Since $\mu_{2m+1}$ is the transitive closure of
	$\mu_{2m}$, we have that $(a,c)\in \mu_{2m+1}\subseteq \mu$.
	Hence $\mu$ is an equivalence relation. 
	
	Let $(a,b)\in \mu$ and $c,d\in M$. There exists a positive integer $k$ such that $(a,b)\in \mu_{4k+3}$. Thus, $(c+a+d,c+b+d)\in \mu_{4(k+1)}\subseteq \mu$. Hence, $\mu$ is a congruence on $(M,+)$. 
	
	Let $(c,c')\in\mu$ and $a,b\in M$ be such that $(c+a,c'+b)\in \mu$. Since $\mu$ is a congruence on $(M,+)$, we have that 
	$(c'+a,c+a)\in \mu$. Hence, $(c'+a,c'+b)\in \mu$. There exists a positive integer $m$ such that $(c'+a,c'+b)\in \mu_{4m+3}$. Hence $(a,b)\in \mu_{4(m+1)}\subseteq \mu$. Therefore, $(M,+)/\mu$ is a left cancellative monoid.
	
	Let $(a,b)\in \mu$ and $c,d\in M$. There exists a positive integer $k$ such that $(a,b)\in \mu_{4k+1}$. It follows that $(\lambda'_{d}(a\circ c),\lambda'_d(b\circ c))\in \mu_{4k+2}$ and $(d\circ a\circ c, d\circ b\circ c)=(d+\lambda'_{d}(a\circ c),d+\lambda'_d(b\circ c))\in \mu_{4(k+1)}\subseteq \mu$. Hence, $\mu$ is a congruence on $(M,\circ)$.
	
	Let $(c,c')\in\mu$ and $a,b\in M$ be such that $(c\circ a,c'\circ b)\in \mu$. Since $\mu$ is a congruence on $(M,\circ)$, we have that 
	$(c'\circ a,c\circ a)\in \mu$. Hence $(c'+\lambda_{c'}(a),c'+\lambda_{c'}(b))=(c'\circ a,c'\circ b)\in \mu$. Since $(M,+)/\mu$ is a left cancellative monoid we get that $(\lambda_{c'}(a),\lambda_{c'}(b))\in \mu$. Now there exists a positive integer $m$ such that
	$(\lambda_{c'}(a),\lambda_{c'}(b))\in \mu_{4m+1}$, and thus $(a,b)\in \mu_{4m+2}\subseteq \mu$. Therefore $(M,\circ)/\mu$ is a left cancellative monoid.

	Let $(a,b),(c,d)\in \mu$. Since $\mu$ is a congruence on $(M,\circ)$, we have that 
	$$(c+\lambda'_c(x),d+\lambda'_d(x))=(c\circ x, d\circ x)\in \mu,$$
	for all $x\in M$.
	Since $(M,+)/\mu$ is a left cancellative monoid, we get that $(\lambda'_c(x),\lambda'_d(x))\in \mu$, for all $x\in M$. For $x=(\lambda'_c)^{-1}(y)$, we have that
	$$(y,\lambda'_d(\lambda'_c)^{-1}(y))\in \mu,$$ 
	for all $y\in M$. Thus, there exists a positive integer $m$ such that $(y,\lambda'_d(\lambda'_c)^{-1}(y))\in \mu_{4m+1}$. Hence
	$((\lambda'_d)^{-1}(y),(\lambda'_c)^{-1}(y))\in \mu_{4m+2}$. Therefore,
	$$((\lambda'_d)^{-1}(y),(\lambda'_c)^{-1}(y))\in \mu,$$ 
	for all $y\in M$. Now there exists a positive integer $k$ such that
	$$((\lambda'_d)^{-1}(a),(\lambda'_c)^{-1}(a)), (\lambda'_d(a),\lambda'_c(a)), (a,b)\in \mu_{4k+1}.$$
	Hence,
	$$((\lambda'_c)^{-1}(a),(\lambda'_d)^{-1}(a)), ((\lambda'_d)^{-1}(a),(\lambda'_d)^{-1}(b)), (\lambda'_c(a),\lambda'_d(a)) (\lambda'_d(a),\lambda'_d(b))\in \mu_{4k+2},$$
	and thus,
	$$(\lambda'_c(a),\lambda'_d(b)),((\lambda'_c)^{-1}(a),(\lambda'_d)^{-1}(b))\in \mu_{4k+3}.$$
	Therefore,
	$$(\lambda'_c(a),\lambda'_d(b)),((\lambda'_c)^{-1}(a),(\lambda'_d)^{-1}(b))\in \mu,$$
	for all $(a,b),(c,d)\in \mu$, and the result follows.
\end{proof}

With the assumptions and notations as in  Example~\ref{ex1},  let $\bar
M=M/\mu$ and let $ M\longrightarrow \bar M: a\mapsto \overline{a}$ be
the natural  map. Let $\bar\lambda\colon (\bar
M,\circ)\longrightarrow \Aut(\bar M,+)$ be the map defined by
$\bar\lambda(\bar a)=\bar\lambda_{\bar a}$ and $\bar\lambda_{\bar
	a}(\bar b)=\overline{\lambda'_a(b)}$, for all $a,b\in M$.

Note that $\bar\lambda$ is well-defined, because if $\bar c=\bar a$
and $\bar d=\bar b$, then, by Lemma~\ref{congruence},
$$\overline{\lambda'_a(b)}=\overline{\lambda'_c(d)}.$$
Now it is easy to check that $\bar\lambda_{\bar a}\in \Aut(\bar
M,+)$, in fact $(\bar\lambda_{\bar a})^{-1}\colon \bar M\longrightarrow \bar M$ is the map defined by $(\bar\lambda_{\bar a})^{-1}(\bar b)=\overline{(\lambda'_a)^{-1}(b)}$, which also  is well-defined by Lemma~\ref{congruence}. Furthermore, by Lemma~\ref{congruence}, $(\bar M,\circ)$ is left cancellative and $\bar\lambda$ is a homomorphism such that $\bar
a\circ\bar b=\bar a+\bar\lambda_{\bar a}(\bar b)$, for all $a,b\in
M$. 

Let $\bar \phi\colon \bar M\times\bar M\longrightarrow \bar M$ be
the map defined by $\bar\phi(\bar a,\bar b)=\bar\lambda_{\bar
	a}(\bar b)$, for all $a,b\in M$. Then $(\bar M,+,\circ,\bar \phi)$
is a left semi-truss.

By \cite[Lemma~5.8]{CJV}, the left cancellative monoid $(\bar M,+)$ satisfies that
for all $\bar a,\bar b\in\bar M$   there exists a unique $\bar c\in
	\bar M$ (denoted as $c(\bar{a},\bar{b})$) such that $\bar a+\bar b=\bar b+\bar c$.  
	Hence, from \cite[Proposition 5.4]{CJV}, we have the following corollary.

\begin{corollary}  Let $(X,r)$ be a left non-degenerate set-theoretic solution of the YBE. Let $\mu$ be the congruence on  $M=(M(X,r'),+)$ defined above.
	Then $(\bar M,+,\circ,\bar \phi)$ is a  left semi-truss   with $\bar M +\bar a \subseteq\bar a +\bar M$ for all $\bar a \in \bar M$ and 
	with $\bar\phi(\bar a,\bar b)=\bar\lambda_{\bar
		a}(\bar b)$, for all $\bar a,\bar b\in \bar M$.
	Furthermore,  $(\bar M,\bar r)$, where
	$$\bar r(\bar a,\bar b)=(\bar \lambda_{\bar a}(\bar b),\bar \lambda^{-1}_{\bar \lambda_{\bar a}(\bar b)}(c(\bar a,\bar \lambda_{\bar a}(\bar b)))),$$
	for all $\bar a,\bar b\in \bar M$, is a left non-degenerate
	set-theoretic solution of the YBE.
	In particular, $(\bar X, \bar r_{| \bar X^2})$ is a left non-degenerate solution on the image $\bar X$ of $X$ in $\bar M$.
\end{corollary}

\section{Addendum}
In this section, we will generalize the first part of \cite[Proposition 4.2]{JKV}.  
Let $\eta$ be the left cancellative congruence on $(M,+)$, defined in \cite{CJV}. 
For a left non-degenerate solution $(X,r)$, we will define the (least)  left cancellative congruence on $(M,\circ)$, say $\nu$, and show that $\eta=\nu$ and $\lambda'_a=\lambda'_b$, for all $(a,b)\in \eta$, in case the solution is bijective and (left and right) non-degenerate. 
We again will follow the notation of \cite{CJV}.

Let $\nu$ be the left cancellative congruence on $(M,\circ)$, that is,
$\nu$ is the smallest congruence such that $\bar M=(M,\circ)/\nu$ is a
left cancellative monoid.

We shall give  a    description of the elements in $\nu$. Let
$$\nu_0=\{ (a,b)\in M^2\mid \exists c\in M\mbox{ such that }c\circ a=c\circ b\}.$$
Note that $\nu_0$ is a reflexive and symmetric binary relation on
$M$. Let $\nu_1$ be its transitive closure, that is
$$\nu_1=\{(a,b)\in M^2\mid \exists a_1,\dots,a_n\in M \mbox{ such that }(a,a_1),(a_1,a_2),\dots ,(a_n,b)\in \nu_0\}.$$
Thus, $\nu_1$ is an equivalence relation on $M$. Let
\begin{align*}\nu_2=&\{ (c\circ a,c\circ b)\in M^2\mid c\in M\mbox{  and
	}(a,b)\in \nu_1\}\\
	&\cup \{ (a,b)\in M^2\mid \exists c\in M\mbox{ such that
	}(c\circ a,c\circ b)\in \nu_1\},
\end{align*}
and for every $m\geq 1$ we
define
$$\nu_{2m+1}=\{(a,b)\in M^2\mid \exists a_1,\dots,a_n\in M \mbox{ such that }(a,a_1),(a_1,a_2),\dots ,(a_n,b)\in \nu_{2m}\}$$
and
\begin{align*}\nu_{2m+2}=&\{ (c\circ a,c\circ b)\in M^2\mid c\in M\mbox{
		 and }(a,b)\in \nu_{2m+1}\}\\
	&\cup \{ (a,b)\in M^2\mid \exists c\in M\mbox{ such that
	}(c\circ a,c\circ b)\in \nu_{2m+1}\},\end{align*}
Note that
$\nu_n\subseteq\nu_{n+1}\subseteq \nu$ for all $n\geq 0$. Let
$\nu'=\cup_{n=0}^{\infty}\nu_n$.

\begin{lemma}\label{lemma1}
	With the above notation we have that $\nu'=\nu$ and
	$\lambda'_a=\lambda'_b$, for all $(a,b)\in \nu$. Furthermore, if $r$ is bijective, then for
	all $z\in M$,
	$$\nu\supseteq \{ ( (\lambda'_z)^{-1}(a),(\lambda'_z)^{-1}(b))\mid (a,b)\in\nu\},$$
	and $\nu$ also is a congruence on $(M,+)$.
\end{lemma}
\begin{proof}
	First we shall prove that $\nu'$ is a congruence on $(M,\circ)$.
	Clearly $\nu'$ is reflexive and symmetric because so is each
	$\nu_n$. Let $a,b,c\in M$ be such that $(a,b),(b,c)\in \nu'$. There
	exists a positive integer $m$ such that $(a,b),(b,c)\in
	\nu_{2m}$. Since $\nu_{2m+1}$ is the transitive closure of
	$\nu_{2m}$, we have that $(a,c)\in \nu_{2m+1}\subseteq \nu'$.
	Hence $\nu'$ is an equivalence relation. Note that every $\nu_n$
	satisfies that $(x\circ z,y\circ z)\in\nu_n$,  for all $(x,y)\in \nu_n$.
	Thus $(a\circ c,b\circ c)\in \nu_{2m}\subseteq \nu'$. Since $(a,b)\in
	\nu_{2m}\subseteq\nu_{2m+1}$, we have that $(c\circ a,c\circ b)\in \nu_{2m+2}\subseteq
	\nu'$. Therefore, $\nu'$ is a congruence.
	
	Let $a,b,c,c'\in M$ be elements such that $(c,c'),(c\circ a,c'\circ b)\in \nu'$. Since $\nu'$ is a congruence on $(M,\circ)$, $(c'\circ b,c\circ b)\in \nu'$. Hence $(c\circ a,c\circ b)\in\nu'$. There exists a positive integer $t$ such that $(c\circ a,c\circ b)\in\nu_{2t+1}.$
	Thus $(a,b)\in \nu_{2t+2}\subseteq \nu'$. Hence $(M,\circ)/\nu'$ is a
	left cancellative monoid. Since $\nu'\subseteq \nu$, we have
	$\nu'=\nu$ by the definition of $\nu$.
	
	Let $(a,b)\in \nu_0$. Then there exists $c\in M$ such that
	$c\circ a=c\circ b$. Hence,
	$$\lambda'_{c}\lambda'_{a}=\lambda'_{c}\lambda'_{b}$$
	and thus $\lambda'_a=\lambda'_b$, for all $(a,b)\in
	\nu_0$. Let $n>0$ and suppose that $\lambda'_a=\lambda'_b$, for all $(a,b)\in
	\nu_{n-1}$. If $n-1$ is even, then for every $(a,b)\in \nu_n$ there exist $(a,c_1), (c_1,c_2),\dots, (c_k,b)\in \nu_{n-1}$. By the induction hypothesis $\lambda'_a=\lambda'_{c_1}=\cdots =\lambda'_{c_k}=\lambda'_b$. If $n-1$ is odd and $(a,b)\in \nu_n$, then either $(a,b)=(c\circ a',c\circ b')$, for some $c\in M$ and $(a',b')\in \nu_{n-1}$, or there exists $c\in M$ such that $(c\circ a,c\circ b)\in\nu_{n-1}$. In the first case, by the induction hypothesis, we have that
	$$\lambda'_a=\lambda'_{c\circ a'}=\lambda'_c\lambda'_{a'}=\lambda'_c\lambda'_{b'}=\lambda'_{c\circ b'}=\lambda'_b.$$
	In the second case, by the induction hypothesis, we have that  
	$$\lambda'_c\lambda'_{a}=\lambda'_{c\circ a}=\lambda'_{c\circ b}=\lambda'_c\lambda'_{b},$$
	and thus $\lambda'_a=\lambda'_b$. Hence, we get that $\lambda'_a=\lambda'_b$, for all $(a,b)\in \nu_n$. Hence, by induction, we have that
	$\lambda'_a=\lambda'_b$, for all $(a,b)\in \nu$.
	
	Suppose that $r$ is bijective. By Example \ref{ex1}, we have that $M+a=a+M$, for all $a\in M$. 
	Let $(a,b)\in \nu_0$. Then there exists $c\in M$ such that
	$c\circ a=c\circ b$. Let $y\in M$. We have that there exists $z\in M$ such that $z+c=c+y$. Hence,
	\begin{align*}(\lambda'_z)^{-1}(c\circ a)=&(\lambda'_z)^{-1}(c+\lambda'_c(a))\\
		=&(\lambda'_z)^{-1}(c)+(\lambda'_z)^{-1}\lambda'_c(a)\\
		=&(\lambda'_z)^{-1}(c)\circ (\lambda'_{(\lambda'_z)^{-1}(c)})^{-1}(\lambda'_z)^{-1}\lambda'_c(a)\\
		=&(\lambda'_z)^{-1}(c)\circ (\lambda'_{z\circ (\lambda'_z)^{-1}(c)})^{-1}\lambda'_c(a)\\
		=&(\lambda'_z)^{-1}(c)\circ (\lambda'_{z+c})^{-1}\lambda'_c(a)\\
		=&(\lambda'_z)^{-1}(c)\circ (\lambda'_{c+y})^{-1}\lambda'_c(a)\\
		=&(\lambda'_z)^{-1}(c)\circ (\lambda'_{c\circ(\lambda'_c)^{-1}(y)})^{-1}\lambda'_c(a)\\
		=&(\lambda'_z)^{-1}(c)\circ (\lambda'_{(\lambda'_c)^{-1}(y)})^{-1}(\lambda'_c)^{-1}\lambda'_c(a)\\
		=&(\lambda'_z)^{-1}(c)\circ (\lambda'_{(\lambda'_c)^{-1}(y)})^{-1}(a).
	\end{align*}
	Since $c\circ a=c\circ b$, we have that
	$$(\lambda'_z)^{-1}(c)\circ (\lambda'_{(\lambda'_c)^{-1}(y)})^{-1}(a)=(\lambda'_z)^{-1}(c)\circ (\lambda'_{(\lambda'_c)^{-1}(y)})^{-1}(b).
	$$
We get that	
		$$((\lambda'_{(\lambda'_c)^{-1}(y)})^{-1}(a),(\lambda'_{(\lambda'_c)^{-1}(y)})^{-1}(b))\in \nu_0,
		$$
		for all $y\in M$. Hence,
		$$\nu_0\supseteq \{ ( (\lambda'_z)^{-1}(a),(\lambda'_z)^{-1}(b))\mid (a,b)\in\nu_0\},$$
		for	all $z\in M$.
Let $n$ be a positive integer and suppose that 	
	$$\nu_{n-1}\supseteq \{ ( (\lambda'_z)^{-1}(a),(\lambda'_z)^{-1}(b))\mid (a,b)\in\nu_{n-1}\},$$
for	all $z\in M$. Let $(a,b)\in \nu_n$. If $n$ is odd, then there exist $(a,c_1),(c_1,c_2),\dots, (c_k,b)\in\nu_{n-1}$. By the induction hypothesis,
$$( (\lambda'_z)^{-1}(a),(\lambda'_z)^{-1}(c_1)),( (\lambda'_z)^{-1}(c_1),(\lambda'_z)^{-1}(c_2)),\dots ,( (\lambda'_z)^{-1}(c_k),(\lambda'_z)^{-1}(b))\in \nu_{n-1}.$$
Hence $( (\lambda'_z)^{-1}(a),(\lambda'_z)^{-1}(b))\in \nu_n$, in this case. If $n$ is even, then either $(a,b)=(c\circ a',c\circ b')$, for some $c\in M$ and $(a',b')\in \nu_{n-1}$, or there exists $c\in M$ such that 
$(c\circ a,c\circ b)\in\nu_{n-1}$. In the first case,
$$(\lambda'_z)^{-1}(a)=(\lambda'_z)^{-1}(c)\circ (\lambda'_{(\lambda'_c)^{-1}(y)})^{-1}(a'),$$
and
$$(\lambda'_z)^{-1}(b)=(\lambda'_z)^{-1}(c)\circ (\lambda'_{(\lambda'_c)^{-1}(y)})^{-1}(b'),$$
where $z+c=c+y$. Hence, by the induction hypothesis,  $( (\lambda'_z)^{-1}(a),(\lambda'_z)^{-1}(b))\in \nu_n$, in this case. In the second case, by the induction hypothesis,
$$( (\lambda'_z)^{-1}(c\circ a),(\lambda'_z)^{-1}(c\circ b))\in \nu_{n-1}.$$
Since $(\lambda'_z)^{-1}(c\circ a)=(\lambda'_z)^{-1}(c)\circ (\lambda'_{(\lambda'_c)^{-1}(y)})^{-1}(a)$, we have that
$$((\lambda'_{(\lambda'_c)^{-1}(y)})^{-1}(a),(\lambda'_{(\lambda'_c)^{-1}(y)})^{-1}(b))\in \nu_n.$$
Since $M+c=c+M$,
$$((\lambda'_{(\lambda'_c)^{-1}(y)})^{-1}(a),(\lambda'_{(\lambda'_c)^{-1}(y)})^{-1}(b))\in \nu_n,$$
for all $y\in M$. Hence,
$$\nu_{n}\supseteq \{ ( (\lambda'_z)^{-1}(a),(\lambda'_z)^{-1}(b))\mid (a,b)\in\nu_{n}\},$$
for	all $z\in M$. By induction, we get that
$$\nu\supseteq \{ ( (\lambda'_z)^{-1}(a),(\lambda'_z)^{-1}(b))\mid (a,b)\in\nu\},$$
for	all $z\in M$.

	Let $(a,b)\in \nu$. Then for every $c\in M$, we have that
	$$(c+a,c+b)=(c\circ (\lambda'_c)^{-1}(a),c\circ (\lambda'_c)^{-1}(b))\in \nu.$$
	Since $\lambda'_a=\lambda'_b$,  we have that
	$$(a+c,b+c)=(a\circ(\lambda'_a)^{-1}(c),b\circ(\lambda'_b)^{-1}(c))=(a\circ(\lambda'_a)^{-1}(c),b\circ(\lambda'_a)^{-1}(c))\in \nu.$$
	Hence $\nu$ is a congruence on $(M,+)$, and the result follows.
\end{proof}

 In order to prove the main result of this section, we first show that, for left non-degenerate set-theoretic solutions of the YBE, the maps $\lambda$ and $\lambda'$ are equal.
Here $\lambda$ is the unique monoid homomorphism  
$M\rightarrow \Map (M,M): a\mapsto \lambda_a$ defined 
in \cite[Theorem 2.1]{CJV} such that 
$\lambda_b (a\circ c) =\lambda_b(a) \circ \lambda_{\rho_a(b)}(c)$ 
and $\rho_{b}(c\circ a) = \rho_{\lambda_a(b)}(c) \circ \rho_b (a)$, 
where also $\rho : M \rightarrow \Map (M,M)$ is the monoid anti-homomorphism 
defined in \cite[Theorem 2.1]{CJV}.
This result comes from \cite{CJKVV}, but for completeness sake we include a proof.

\begin{lemma}\label{lemma2}
	Let $(X,r)$ be a set-theoretic solution of the YBE. Let $M=M(X,r)$ and $M'=A(X,r)$. As usual, write  $r(x,y)=(\sigma_x(y),\gamma_y(x))$. Then, $\lambda'_a(\pi(b))=\pi(\lambda_a(b))$, for all $a,b\in M$, where  $\pi\colon M\rightarrow M'$ is the unique $1$-cocycle with respect the left action $\lambda'$ such that $\pi(x)=x$, for all $x\in X$. Furthermore, if $(X,r)$ is left non-degenerate, then, with the identification of $M$ and $M'$ in Example \ref{ex1}, $\lambda'_a(b)=\lambda_a(b)$, for all $a,b\in M$. In particular,
	\begin{equation}\label{gamma}\lambda'_x(x_1\circ\cdots\circ x_k\circ a)=\lambda'_x(x_1\circ\cdots\circ x_k)\circ\lambda'_{\gamma_{x_k}\cdots\gamma_{x_1}(x)}(a),
	\end{equation}
	for all $x,x_1,\dots, x_k\in X$ and $a\in M$.
\end{lemma}

\begin{proof} The existence and uniqueness of $\pi$ is proven in \cite[Proposition 3.2]{CJV}.
Let $b\in M$. There exist a non-negative integer $k$ and $x_1,\dots ,x_k\in X$ such that $b=x_1\circ \cdots\circ x_k$. We first prove that
$\lambda'_x(\pi(b))=\pi(\lambda_x(b))$, for all $x\in X$, by induction on $k$. If $k=0$, then $\pi(1)=0$ and by the definition of $\lambda$, $\lambda_x(1)=1$. Hence $\lambda'_x(\pi(1))=\lambda'_x(0)=0=\pi(1)=\pi(\lambda_x(1))$. For $k=1$, \[\pi(\lambda_x(x_1))=\sigma_x(x_1)=\sigma_x(\pi(x_1))=\lambda'_x(\pi(x_1)).\]
Suppose that $k>1$ and we have proven the result for words in $M(X,r)$ of length at most $k-1$. By the definition of $\lambda$, \cite[Theorem 2.1]{CJV} and the induction hypothesis, we have
\begin{align*}
	\pi(\lambda_x(b)) &=\pi(\lambda_x(x_1\circ \cdots \circ x_k))\\
	& =\pi(\lambda_x(x_1)\circ\lambda_{\rho_{x_1}(x)}(x_2\circ\cdots\circ x_k))\\
	& =\pi(\lambda_x(x_1))+\lambda'_{\lambda_x(x_1)}(\pi(\lambda_{\rho_{x_1}(x)}(x_2\circ\cdots \circ x_k)))\\
	& =\lambda'_x(\pi(x_1))+\lambda'_{\lambda_x(x_1)}(\lambda'_{\rho_{x_1}(x)}(\pi(x_2\circ\cdots\circ x_k)))\\
	& =\lambda'_x(\pi(x_1))+\lambda'_x(\lambda'_{x_1}(\pi(x_2\circ\cdots\circ x_k)))\\
	& =\lambda'_x(\pi(x_1)+\lambda'_{x_1}(\pi(x_2\circ\cdots\circ x_k)))\\
	& =\lambda'_x(\pi(x_1\circ\cdots \circ x_k))\\
	& =\lambda'_x(\pi(b)).
\end{align*}
Hence, by induction $\lambda'_x(\pi(b))=\pi(\lambda_x(b))$, for all $x\in X$ and $b\in M$.
Using that both $\lambda$ and $\lambda'$ are homomorphisms, we obtain $\lambda'_{a}(\pi(b))=\pi(\lambda_{a}(b))$
for all $a,b\in M$.	
	
Suppose that $(X,r)$ is left non-degenerate. Then with the identification of $M$ and $M'$ in Example \ref{ex1}, we have that 
$\lambda'_a(b)=\lambda_a(b)$, for all $a,b\in M$. In this case, by \cite[theorem 2.1]{CJV},
\begin{align*}\lambda'_x(x_1\circ\cdots\circ x_k\circ a)=&\lambda_x(x_1\circ\cdots\circ x_k\circ a)\\
	=&\lambda_x(x_1\circ \cdots \circ x_k)\circ\lambda_{\rho_{x_1\circ\cdots\circ x_k}(x)}(a)\\
	=&\lambda_x(x_1\circ \cdots \circ x_k)\circ\lambda_{\gamma_{x_k}\cdots \gamma_{x_1}(x)}(a)\\
	=&\lambda'_x(x_1\circ x_2\circ\cdots\circ x_k)\circ \lambda'_{\gamma_{x_k}\cdots \gamma_{x_1}(x)}(a),
\end{align*}
for all $x,x_1,\dots ,x_k\in X$ and $a\in M$. Hence, (\ref{gamma}) follows. 
\end{proof}

\begin{proposition}\label{q}
	Let $(X,r)$ be a bijective (left and right) non-degenerate set-theoretic solution of the YBE. Let $M=M(X,r)$. As usual, write  $r(x,y)=(\sigma_x(y),\gamma_y(x))$.  Let $\nu$ be the left cancellative congruence on $(M,\circ)$, and let $\eta$ be the left cancellative congruence on $(M,+)$. Then $\eta=\nu$ and thus, for every $z\in M$,
	$$\nu=\{ (\lambda'_z(a),\lambda'_z(b))\mid (a,b)\in \nu\}=\{((\lambda'_z)^{-1}(a),(\lambda'_z)^{-1}(b))\mid (a,b)\in \nu\}.$$
	Furthermore $\lambda'_a=\lambda'_b$, for all $(a,b)\in \eta$.
	\end{proposition}

\begin{proof}
	From the proof of Lemma \ref{lemma1}, we know that for
	all $z\in M$,
	$$\nu_0\supseteq \{ ( (\lambda'_z)^{-1}(a),(\lambda'_z)^{-1}(b))\mid (a,b)\in\nu_0\},$$
	and $\nu$ also is a congruence on $(M,+)$.

	Let $(a,b)\in \nu_0$. Then there exists $c\in M$ such that
	$c\circ a=c\circ b$. There exist $x_1,\dots ,x_k\in X$ such that $c=x_1\circ\cdots\circ x_k$.  Let $x\in X$. By (\ref{gamma})  (in Lemma \ref{lemma2}), we have that
	$$\lambda'_x(c\circ a)=\lambda'_x(c)\circ\lambda'_{\gamma_{x_k}\cdots \gamma_{x_1}(x)}(a).$$
	Hence
	$$\lambda'_x(c)\circ\lambda'_{\gamma_{x_k}\cdots \gamma_{x_1}(x)}(a)=\lambda'_x(c)\circ\lambda'_{\gamma_{x_k}\cdots \gamma_{x_1}(x)}(b),$$
	for all $x\in X$.
Hence, $(\lambda'_{\gamma_{x_k}\cdots \gamma_{x_1}(x)}(a), \lambda'_{\gamma_{x_k}\cdots \gamma_{x_1}(x)}(b)) \in \nu_0$, for all $x \in X$. 	
Since $(X,r)$ is right non-degenerate, and thus all $\gamma_{x_i}$ are bijective, we obtain that 
$(\lambda'_y(a),\lambda'_y(b))\in\nu_0$, for all $y\in X$. Therefore
$(\lambda'_z(a),\lambda'_z(b))\in\nu_0$, for all $z\in M$. Since $((\lambda'_z)^{-1}(a),(\lambda'_z)^{-1}(b))\in\nu_0$, for all $z\in M$,
we get that
$$\nu_0=\{ (\lambda'_z(a),\lambda'_z(b))\mid (a,b)\in \nu_0\}=\{((\lambda'_z)^{-1}(a),(\lambda'_z)^{-1}(b))\mid (a,b)\in \nu_0\},$$
for all $z\in M$. We shall prove by induction on $n$ that
\begin{equation}\label{induction}
	\nu_n=\{ (\lambda'_z(a),\lambda'_z(b))\mid (a,b)\in \nu_n\}=\{((\lambda'_z)^{-1}(a),(\lambda'_z)^{-1}(b))\mid (a,b)\in \nu_n\},
\end{equation}
	for all $z\in M$ and all non-negative integers $n$. 	Suppose that $n>0$  and (\ref{induction}) is true for $n-1$. From the proof of Lemma \ref{lemma1}, we know that for
	all $z\in M$,
	$$\nu_n\supseteq \{ ( (\lambda'_z)^{-1}(a),(\lambda'_z)^{-1}(b))\mid (a,b)\in\nu_n\}.$$
	Let $(a,b)\in \nu_n$, $z\in M$. If $n$ is odd, then there exist $(a,c_1),(c_1,c_2),\dots ,(c_k,b)\in \nu_{n-1}$. Hence,
	$$(\lambda'_z(a),\lambda'_z(c_1)), (\lambda'_z(c_1),\lambda'_z(c_2)),\dots, (\lambda'_z(c_k),\lambda'_z(b))\in\nu_{n-1},$$
	and thus
	$(\lambda'_z(a),\lambda'_z(b))\in \nu_n$, in this case.
	If $n$ is even then either $(a,b)=(c\circ a',c\circ b')$, for some $c\in M$ and $(a',b')\in \nu_{n-1}$, or there exists $c\in M$ such that $(c\circ a,c\circ b)\in\nu_{n-1}$. Put $c=x_1 \circ \cdots \circ x_k$.
 	In the first case,
    by the previous lemma, we get, $(\lambda'_z(a),\lambda'_z(b))= (\lambda'_z(c\circ a'),\lambda'_z(c\circ b'))\stackrel{\eqref{gamma}}{=} (\lambda'_z(c)\circ\lambda'_{\gamma_{x_k} \cdots \gamma_{x_1}(z)}(a'),\lambda'_z(c)\circ\lambda'_{\gamma_{x_k} \cdots \gamma_{x_1}(z)}(b'))$. 
	By the induction hypothesis, and since $(a',b') \in \nu_{n-1}$, also $(\lambda'_{\gamma_{x_k} \cdots \gamma_{x_1}(z)}(a'),\lambda'_{\gamma_{x_k} \cdots \gamma_{x_1}(z)}(b')) \in \nu_{n-1}$, and then $(\lambda'_z(a),\lambda'_z(b))=(\lambda'_z(c)\circ\lambda'_{\gamma_{x_k} \cdots \gamma_{x_1}(z)}(a'),\lambda'_z(c)\circ\lambda'_{\gamma_{x_k} \cdots \gamma_{x_1}(z)}(b'))\in \nu_n$ (because $n$ is even).  
	In the second case, 
	by (\ref{gamma}),
	$$\lambda'_x(c\circ a)=\lambda'_x(x_1\circ\cdots\circ x_k\circ a)=\lambda'_x(x_1\circ\cdots\circ x_k)\circ\lambda'_{\gamma_{x_k}\cdots\gamma_{x_1}(x)}(a),$$ 
	for all $x\in X$. By the induction hypothesis,
	$$(\lambda'_x(c\circ a),\lambda'_x(c\circ b))\in\nu_{n-1}, $$
	for all $x\in X$. Hence,
	$$(\lambda'_{\gamma_{x_k}\cdots\gamma_{x_1}(x)}(a),\lambda'_{\gamma_{x_k}\cdots\gamma_{x_1}(x)}(b))\in \nu_n,$$
	for all $x\in X$. Since $(X,r)$ is right non-degenerate, we have that
		$$(\lambda'_{y}(a),\lambda'_{y}(b))\in \nu_n,$$
	for all $y\in X$. Hence, $(\lambda'_{z}(a),\lambda'_{z}(b))\in \nu_n$, for all $z\in M$. Since
	$$\nu_n\supseteq \{ ( (\lambda'_z)^{-1}(a),(\lambda'_z)^{-1}(b))\mid (a,b)\in\nu_n\},$$
	we get that  
	$$\nu_n=\{ ( \lambda'_z(a),\lambda'_z(b))\mid (a,b)\in\nu_n\}= \{ ( (\lambda'_z)^{-1}(a),(\lambda'_z)^{-1}(b))\mid (a,b)\in\nu_n\},$$
for all $z\in M$. Hence, by induction, 
	$$\nu=\{ ( \lambda'_z(a),\lambda'_z(b))\mid (a,b)\in\nu\}= \{ ( (\lambda'_z)^{-1}(a),(\lambda'_z)^{-1}(b))\mid (a,b)\in\nu\},$$
for all $z\in M$.

Let $a,b,c,c'\in M$  be such that $(c,c'),(c+a,c'+b)\in \nu$. Since $\nu$ is a congruence on $(M,+)$, $(c'+b,c+b)\in\nu$. Hence $(c+a,c+b)\in \nu$. Then,
$(c\circ (\lambda'_c)^{-1}(a),c\circ (\lambda'_c)^{-1}(b))=(c+a,c+b)\in\nu$.  Hence, $( (\lambda'_c)^{-1}(a), (\lambda'_c)^{-1}(b))\in\nu$ and thus $(a,b)\in \nu$. Therefore, $(M,+)/\nu$ is left cancellative  and thus clearly  $\eta\subseteq \nu$.

By Lemma \ref{lemma1}, $\lambda'_a=\lambda'_b$, for all $(a,b)\in\eta \subseteq \nu$. Let $(a,b)\in \eta$ and let $c\in M$. By Lemma \ref{eta}, we have that
$(c\circ a,c\circ b)=(c+\lambda'_c(a),c+\lambda'_c(b))\in\eta$, and since $\lambda'_a=\lambda'_b$, we have that
$(a\circ c, b\circ c)=(a+\lambda'_a(c),b+\lambda'_b(c))=(a+\lambda'_a(c),b+\lambda'_a(c))\in\eta$. Hence, $\eta$ is a congruence on $(M,\circ)$.   Let $a,b,c,c'\in M$ be such that $(c,c'),(c\circ a,c'\circ b)\in \eta$. Then, $(c+\lambda'_c(a), c'+\lambda'_{c'}(b))=(c\circ a,c'\circ b)\in\eta$. Since $\lambda'_c=\lambda'_{c'}$, we have that
$$(c+\lambda'_c(a), c'+\lambda'_{c}(b)), (c'+\lambda'_c(b), c+\lambda'_{c}(b))\in\eta$$	
and then $(c+\lambda'_c(a), c+\lambda'_{c}(b))\in\eta$.
Hence,
$(\lambda'_c(a),\lambda'_c(b))\in \eta$. By Lemma \ref{eta}, $(a,b)\in \eta$. Therefore, $(M,\circ)/\eta$ is left cancellative and
$\nu\subseteq\eta$. So, $\eta=\nu$ and the result follows. 	
\end{proof}

\vspace{30pt}
 \noindent \begin{tabular}{llllllll}
  Ferran Ced\'o && Eric Jespers\\
 Departament de Matem\`atiques &&  Department of Mathematics \\
 Universitat Aut\`onoma de Barcelona &&  Vrije Universiteit Brussel \\
08193 Bellaterra (Barcelona), Spain    && Pleinlaan 2, 1050 Brussel, Belgium \\
 cedo@mat.uab.cat && Eric.Jespers@vub.be \\ \\
 Charlotte Verwimp && \\ Department of
Mathematics &&
\\  Vrije Universiteit Brussel &&\\
Pleinlaan 2, 1050 Brussel, Belgium  &&\\
Charlotte.Verwimp@vub.be &&
\end{tabular}

\end{document}